\documentclass[11pt]{amsart}
\usepackage{amscd, amssymb}

\newtheorem{Theorem}{Theorem}
\newtheorem{Proposition}{Proposition}

\newtheorem{Def/Thm}{Definition/Theorem}
\newtheorem{Corollary}{Corollary}
\newtheorem{Lemma}{Lemma}
\theoremstyle{remark}
\newtheorem{Remark}{Remark}
\newtheorem{Example}{Example}

\newcommand{\ra }{\rightarrow}

\newcommand{\Spec}{{\mathrm{Spec}}}

\newcommand{\PP }{{\mathbb P}}

\newcommand{\CC }{{\mathbb C}}
\newcommand{\ZZ }{{\mathbb Z}}

\def\beq{\begin{equation}}
\def\eeq{\end{equation}}
\def\bea{\begin{eqnarray}}
\def\eea{\end{eqnarray}}
\def\bef{\begin{figure}}
\def\enf{\end{figure}}

\def\ba{\begin{array}}
\def\ea{\end{array}}
\def\bce{\begin{center}}
\def\ece{\end{center}}

\def\IC{{\relax\hbox{$\inbar\kern-.3em{\rm C}$}}}
\def\ID{\relax{\rm I\kern-.18em D}}
\def\IE{\relax{\rm I\kern-.18em E}}
\def\IF{\relax{\rm I\kern-.18em F}}
\def\IG{\relax\hbox{$\inbar\kern-.3em{\rm G}$}}
\def\IGa{\relax\hbox{${\rm I}\kern-.18em\Gamma$}}
\def\IH{\relax{\rm I\kern-.18em H}}
\def\II{\relax{\rm I\kern-.18em I}}
\def\IK{\relax{\rm I\kern-.18em K}}

\def\IQ{\relax\hbox{$\inbar\kern-.3em{\rm Q}$}}

\begin{document}

\title{rational curves on blowing-ups of projective spaces}

\keywords{rational curves, moduli spaces, rationally connected
varieties, blowing-ups \\
This work was supported by Korea Research Foundation Grant
(KRF-2004-042-C00005). }

\author{Bumsig Kim}
\address{School of Mathematics, Korea Institute for Advanced
Study, 207-43 Cheong nyangni 2-dong, Dongdaemun-gu, Seoul 130-722,
Korea} \email{bumsig@kias.re.kr}

\author{Yongnam Lee}
\address{Department of Mathematics, Sogang University,
Sinsu-dong, Mapo-gu, Seoul 121-742, Korea }
\email{ynlee@sogang.ac.kr}

\author{Kyungho Oh}
\address{Department of Mathematics and Computer Science,
University of Missouri at St. Louis, One University Boulevard, St.
Lous, MO 63121, USA } \email{oh@arch.umsl.edu}

\maketitle

\section{Introduction}
Let $\mathrm{Mor}_\beta (\mathbb{P} ^1, Y)$ denote the moduli
space of morphisms $f$ from a complex projective line $\mathbb{P}
^1$ to a smooth complex projective variety $Y$ such that
$f_*[\mathbb{P} ^1]=\beta$ where $\beta$ is a given second
homology class of $Y$. We study the irreducibility and the
rational connectedness of the moduli space when $Y$ is a
successive blowing-up of a product of projective spaces with a
suitable condition on $\beta$.

To state the Main Theorem proven in this paper, let us introduce
some notation. Let $X=\prod _{k=1}^m \mathbb{P} ^{n_k}$, $X_0=X$
and let $\pi _i: X_i \rightarrow X_{i-1}$, $i=1,...,r$, be a
blowing-up of $X_{i-1}$ along a smooth irreducible subvariety
$Z_i$. Let $E_i^t\subset X_r$ be the total transform $(\pi _i\circ
...\circ \pi _r)^{-1} Z_i$ of the exceptional divisor associated
to $Z_i$ and let $H_k$ be the divisor class coming from the
hyperplane class of the $k$-th projective space $\mathbb{P}
^{n_k}$. Let $m_i=\# \{Z_j \ | \ j< i ,\ \ (\pi _j\circ ...\circ
\pi _r)^{-1}(Z_j)\supset E_i^t\}$. So general points of $Z_i$ are
the $m_i$-th infinitesimal points of $X$. Denote by
$\mathrm{Mor}_\beta (\mathbb{P} ^1, X_r)^\sharp$ the open sublocus
of $\mathrm{Mor}_\beta (\mathbb{P} ^1, X_r)$, consisting of $f$
whose image does not lie on exceptional divisors: $f(\PP
^1)\nsubseteq E_i^t$, $\forall\, i$.

\bigskip

\noindent{\bf Main Theorem.} {\it Assume that $\beta\cdot (\pi ^*
H_k - \sum _{i=1}^r (m_i+1)E_i^t)\ge 0,$ $\forall\, k$; and $\beta
\cdot E_i^t\ge 0,$ $\forall\, i$, where $\pi = \pi _1\circ
...\circ \pi _r$.

\begin{enumerate}

\item The moduli space $\mathrm{Mor}_\beta (\mathbb{P} ^1,
X_r)^\sharp$ consists of free morphisms and is an irreducible
smooth variety of expected dimension.

\item If $Z_i$ are rationally connected for all $i$, then a
projective, birational model of $\mathrm{Mor}_\beta (\mathbb{P}
^1, X_r)^\sharp$ is also rationally connected.

\item The moduli space $\mathrm{Mor}_\beta (\mathbb{P} ^1, X_r)$ is
smooth and $\mathrm{Mor}_\beta (\mathbb{P} ^1, X_r)^\sharp$ is
dense in $\mathrm{Mor}_\beta (\mathbb{P} ^1, X_r)$ if one of the
following conditions hold
\begin{enumerate}
\item All $\pi (E_i^t)$ are points in $X$.

\item All centers $Z_i$ are convex (that is, $H^1(\PP ^1,
g^*T_{Z_i})=0$ for any morphism $g:\PP ^1\ra Z_i$) and $\pi
(E_i^t)$ are disjoint to $\pi (E_j^t)$ for any $i\ne j$.
\end{enumerate}
\end{enumerate}}

\bigskip

Note that the irreducibility (respectively, the rational
connectedness of a projective, birational model) of the morphism
space $\mathrm{Mor}_\beta (\mathbb{P} ^1, X_r)$ implies that of
the moduli space of rational curves $C$ with numerical condition
$[C]=\beta$.

\medskip

Our paper is motivated by two questions:
\begin{enumerate}

 \item If $Y$ is rationally
connected, is $\mathrm{Mor}_\beta (\mathbb{P} ^1, Y)$ also
irreducible and rationally connected? If this does not always
hold, does it hold for special values of $\beta$?

\item For surfaces fibered over $\PP ^1$ with genus 2 fibers (or
more generally hyperelliptic fibers) and fixed numerical
invariants, is the moduli space of such surfaces connected? If
this does not always hold, does it hold for special values of the
numerical invariants?

\end{enumerate}

\medskip

Several authors have studied the first question. The case of
homogeneous spaces was treated by the first author and
Pandharipande \cite{KP}, and by Thomsen \cite{Tom}. The case of
small degree $d$ general hypersurfaces $Y=X_d\subset\mathbb{P} ^n$
was handled by Harris, Roth and Starr \cite{HRS}. The case when
$Y$ is the moduli spaces of rank 2 stable vector bundles, with
fixed determinant of degree 1, on a smooth projective curve of
genus $g\ge 2$ was investigated by Castravet \cite{Cas}. She found
all irreducible components and described the maximal rationally
connected fibration of them.

Let $\overline{M}_{0,n}$ be the moduli space of stable $n$-pointed
rational curve. As a corollary of the Main Theorem, the space
$\mathrm{Mor}_\beta (\mathbb{P} ^1, \overline{M}_{0,n} )$ is
connected for certain values of $\beta$,  since the space
$\overline{M}_{0,n}$ is a successive blowing-up of $(\mathbb{P}
^1)^{n-3}$ along smooth codimension 2 subvarieties (\cite{Keel})
or a successive blowing-up of $\mathbb{P} ^{n-3}$ (\cite{Kap}).
This gives a step toward proving connectedness of the moduli space
of hyperelliptic fibrations over $\PP ^1$ (presumably by replacing
the hyperelliptic fibration by the fibration of quotient by the
hyperelliptic involution, marked by the images of the Weierstrass
points).

\bigskip

When $Y$ is a blowing-up of a product of projective spaces along a
smooth closed (not necessarily irreducible) subvarieties, we prove
a slightly stronger result, Theorem \ref{theoremA} in section 2.
In section 3, we prove the Main Theorem. The key idea of both
proofs is to express the moduli space as a fibration
--- a fiber consists of the morphisms $f$ which pass through given
points of $\coprod _i Z_i$ at given points of domain $\PP ^1$ ---
and then to show that the general fiber is rationally connected
and has the expected dimension under the condition on $\beta$ as
in the Main Theorem (and also as in theorem \ref{theoremA}). When
$Y$ is a successive blowing-up, we will need to utilize Jet spaces
and Jet conditions in order to show the rational connectedness of
the general fiber. Furthermore we apply a result of
Graber-Harris-Starr \cite{GHS}.

\medskip
Throughout the paper, we will employ the well-known results of the
deformation theory of morphisms of curves, as well as the
established notation as in \cite{Kollar}. The complex number field
$\mathbb{C}$ will be the base field.

\bigskip

{\bf Acknowledgements.} Y.L. and K.O. thank staffs in Korea
Institute for Advanced Study for the hospitality during their
visit. The authors are grateful to Igor Dolgachev for raising the
second problem to us. B.K. thanks Jun-Muk Hwang for a useful
comment on free morphisms. Y.L. also thanks A-M. Castravet for
explaining her result. The authors heartily thank the referee for
valuable suggestions and insightful comments.

\section{Blowing-ups along a smooth subvariety}
\subsection{Set-up and a morphism $\sigma$}\label{sigmasubsection}
Let $\pi: \tilde{X}\ra X$ be the blowing-up of a smooth projective
variety $X$ along a smooth closed subvariety $Z$ with the
exceptional divisor $E$.

For a curve class $\beta\in H_2(\tilde{X},\ZZ )$, consider the
evaluation morphism
\[\begin{array}{ccc} \mathrm{ev}: \PP ^1\times \mathrm{Mor}_\beta (\mathbb{P} ^1,
\tilde{X}) &\rightarrow& \tilde{X} \\
(p,\tilde{f}) &\mapsto& \tilde{f}(p). \end{array}\]

Assume that $\beta\cdot E_i\ge 0$ for all $i=1,...,r$ where $E_i$
are the exceptional irreducible divisors over the irreducible
components $Z_i$ of $Z = \coprod _{i=1}^r Z_i$. In general,
$\tilde{f}(\PP ^1)\subset E _i$  does not imply $\tilde{f}(\PP
^1)\cdot E_i <0$. For example, we have:

{\Example\label{ruled-example}
 Consider the blowing-up $\tilde{X}$ of $X=\PP ^3$ along a
curve $Z\cong \PP ^1$ with a normal bundle
$N_{Z/X}=\mathcal{O}(1)\oplus\mathcal{O}(2)$. Then $E = \PP
(\mathcal{O}_Z (1)\oplus \mathcal{O}_{Z}(2))$. If $C$ is a positive
section (resp. the negative section), then by the construction of
the sections and the universal property of the projectivization $E$
of $\mathcal{O}_Z (1)\oplus \mathcal{O}_{Z}(2)$, we see that $C\cdot
E = 1$ (resp. $C\cdot E = 2$).}

\bigskip
This observation forces us to consider   an open subvariety
$\mathrm{Mor}_\beta (\mathbb{P} ^1, \tilde{X})^\sharp$  of
$\mathrm{Mor}_\beta (\mathbb{P} ^1, \tilde{X})$ consisting of
$\tilde{f}$ such that $\tilde{f}(\PP ^1)\nsubseteq E$. Now the
scheme-theoretic intersection $\Gamma _{\pi\circ \mathrm{ev}}\cap
(\PP ^1\times \mathrm{Mor}_\beta (\PP ^1,\tilde{X})^\sharp\times
Z)\cong \mathrm{ev}^{-1}(E)$, where $\Gamma
_{\pi\circ\mathrm{ev}}$ is the graph of the morphism $\pi\circ
\mathrm{ev}$, can be regarded as a closed subscheme of $\PP
^1\times \mathrm{Mor}_\beta (\mathbb{P} ^1, \tilde{X})^\sharp
\times Z$, which is proper and flat over $\mathrm{Mor}_\beta
(\mathbb{P} ^1, \tilde{X})^\sharp $. Thus it induces a natural
morphism
\[\begin{array}{ccc} \sigma: \mathrm{Mor}_\beta (\mathbb{P} ^1,
\tilde{X})^\sharp  &\rightarrow &\prod
_{i=1}^r\mathrm{Hilb}^{e_i}(\mathbb{P} ^1\times
Z _i) \\
\tilde{f} & \mapsto & (\Gamma _{f}\cap (\PP ^1\times
Z_1),...,\Gamma _{f}\cap (\PP ^1\times Z_r)),
\end{array}\] where $f:=\pi\circ \tilde{f}$ and $e_i:=\beta \cdot E_i$.
 Here $\mathrm{Hilb}^0Y$ of a variety $Y$ is defined to be
$\Spec\, \CC$.

\subsection{An exact sequence}
The following lemma shows a sufficient condition of the generic
smoothness of $\sigma$.

{\Lemma\label{exact} Let $\pi: \tilde{X}\ra X$ as above and
suppose that $\tilde{f}(\PP ^1)\nsubseteq E$.
\begin{enumerate} \item There is a natural injective
morphism of sheaves
\[ \tilde{f}^*\pi ^* T_X (-E) \rightarrow
\tilde{f}^*T_{\tilde{X}}.
\]

\item If $\tilde{f}:\PP ^1\rightarrow \tilde{X}$ is transversal to
$E$, then the above injective morphism induces a short exact
sequence
\[ 0 \rightarrow \tilde{f}^*\pi ^* T_X (-E) \rightarrow
\tilde{f}^*T_{\tilde{X}} \rightarrow (\mathrm{id}_{\PP ^1}\times
f)^* T_{\PP ^1\times Z} \rightarrow 0. \] Furthermore, the
associated morphism \[ H^0 (\PP ^1, \tilde{f}^* T_{\tilde{X}}) \ra
(\mathrm{id}_{\PP ^1}\times f)^* T_{\PP ^1\times Z} \] is the
derivative of $\sigma$.

\item For every $\tilde{f}$ satisfying
$h^1 (\PP ^1, \tilde{f}^*\pi ^* T_X (-E))=0$,
the morphism $\sigma$ is smooth at $[\tilde{f}]$.

\end{enumerate}}

\begin{proof}
The first morphism is defined by the pull-back of the extension of
the isomorphism $T_X \cong T_{\tilde{X}}$ away from $E$. To check
the existence of the extension, we consider the blowing-up
$\tilde{X} \rightarrow X$ locally as $(t,{\bf x},{\bf y}) \mapsto
(z_1=t, {\bf z_2}=t{\bf x}, {\bf z_3}={\bf y})$, where $t,\, {\bf
x},\, {\bf y}$ (resp. $z_1$, ${\bf z}_2$, ${\bf z}_3$) is a system
of local parameters of $\tilde{X}$ (resp. $X$), and the bold
letters denote multi-variables. Then the natural morphism of
sheaves
$$\pi ^* T_X (-E) \rightarrow T_{\tilde{X}}$$ defined by
$(t\frac{\partial}{\partial z_1} \mapsto t
\frac{\partial}{\partial t} - \sum x_i \frac{\partial}{\partial
x_i} )$, $(t\frac{\partial}{\partial {\bf z}_2}\mapsto
\frac{\partial}{\partial {\bf x}} )$, and
$(t\frac{\partial}{\partial {\bf z}_3}\mapsto t
\frac{\partial}{\partial {\bf y}} )$ is the extension. The second
morphism is defined by, at $p$ with $f(p)\in Z$,
$$(f_*|_p)^{-1}\oplus\pi _*|_{\tilde{f}(p)}:
T_{\tilde{X}}|_{\tilde{f}(p)} = f_*T_{\PP ^1}|_{p} \oplus
T_E|_{\tilde{f}(p)} \rightarrow T_{\PP ^1}|_p \oplus T_Z
|_{f(p)},$$ where $T_Y|_y$ denotes the tangent space of a variety
$Y$ at a point $y$. Now the rest of the proof is straightforward.
\end{proof}

\bigskip

{\Remark   In fact the proof above shows that there is an exact
sequence
\[ 0 \ra (\mathrm{ev})^*(\pi ^*T_X(-E)) \ra \mathrm{ev}^*
T_{\tilde{X}}\ra (\mathrm{id}_{\PP ^1}\times
\pi\circ\mathrm{ev})^* T_{\PP ^1\times Z}\ra 0 \] over $\PP
^1\times \mathrm{Mor}_\beta (\PP ^1,\tilde{X})^{\circ}$ where
$\mathrm{Mor}_\beta (\PP ^1,\tilde{X})^{\circ}$ is the locus of
all morphisms in $\mathrm{Mor}_\beta (\PP ^1,\tilde{X})$ which are
transversal to $E$.}

{\Remark  In Example \ref{ruled-example}, the 1-1 morphism whose
image is the negative section is not free: The exact sequence
\[0\ra N_{C/E} = \mathcal{O}(-1) \ra N_{C/\tilde{X}}\ra
N_{E/\tilde{X}}|_C= \mathcal{O}(2)\ra 0 \] splits since $C$ is a
section. Therefore (1) of Lemma \ref{exact} is not true in general
without the condition $\tilde{f}(\PP ^1)\nsubseteq E$.}

\bigskip
When $\tilde{f}(\PP ^1)\subseteq E$, instead of Lemma \ref{exact}
we have the following lemma.

 {\Lemma\label{ruled} Assume that $k\le e+1$, $e=E\cdot \tilde{f}_*[\PP ^1]\ge 0$, and
 $\tilde{f}(\PP ^1)\subseteq E$. If $H^1(\PP^1, f^*T_X
(-e-k))=0$ and $H^1(\PP ^1,f^*T_Z (-k))=0$, then we attain $H^1(\PP
^1, \tilde{f}^*T_{\tilde{X}}(-k))=0$.}

\begin{proof}
i) Note that $H^1(\PP ^1,f^*N_{Z/X}(-e-k))=0$ by $H^1(\PP
^1,f^*T_X (-e-k))=0$ and
\[ 0 \ra T_Z \ra T_X |_Z\ra N_{Z/X}\ra 0 .\]

ii) Note that $H^1(\PP ^1,\tilde{f}^*T_\pi (-k))=0$ and $H^1(\PP
^1,\tilde{f}^*T_E (-k))=0$ by i), $H^1(\PP ^1,f^*T_Z(-k))=0$, and
the exact sequences
\[ 0 \ra \mathcal{O} \ra \pi ^*(N_{Z/X})\otimes
\mathcal{O}_{E}(1) \ra T_\pi \ra 0 ;\] \[ 0 \ra T_\pi \ra T_E \ra
\pi ^* T_Z \ra 0 , \] where $T_\pi$ denotes the relative tangent
bundle of $\pi$.

iii) Finally, $H^1(\PP ^1,\tilde{f}^*T_{\tilde{X}}(-k))=0$ by
\[ 0 \ra T_E \ra T_{\tilde{X}}|_E \ra
\mathcal{O}_{\tilde{X}}(E)|_E \ra 0 .\] Here we use the condition
that $-k+e\ge -1$. \end{proof}

\subsection{The fiber of $\sigma$ when $X=\PP
^n$}\label{fibersubsection}

Let $X=\PP ^n$ and denote
$$c_H(\beta )=\beta \cdot (\pi ^*H-E),$$ where $H$ is the
hyperplane class of $\PP ^n$. If $c_H(\beta) \ge -1$, then the
first statement of Lemma \ref{exact} implies the vanishing of the
obstruction $H^1(\PP ^1, \tilde{f}^*T_{\tilde{X}})=0,\ $ if
$\tilde{f}(\PP ^1)\nsubseteq E$, and hence the space
$\mathrm{Mor}_\beta (\PP ^1,\tilde{X})^\sharp$ is smooth. In
addition, the general fiber of a morphism $\sigma$ is smooth and
has the expected dimension due to Lemma \ref{exact}. Here the
expected dimension of the fiber is by definition
\begin{eqnarray*} & & \mathrm{exp.}\,\dim\mathrm{Mor}_\beta (\PP
^1,\tilde{X}) - \dim \prod _{i=1}^r \mathrm{Hilb}^{e_i}(\PP ^1\times
Z_i) \\ &=& \dim \PP H^0(\mathcal{O}_{\PP ^1}(d)\otimes \CC ^{n+1})
-ne = \dim \mathrm{Mor}_{\pi _*\beta} (\PP ^1,X) -ne.\end{eqnarray*}

In the following lemma we investigate the irreducibility of the
fiber of the morphism $\sigma$ for $X=\PP ^n$.

{\Lemma\label{fiber} Suppose that $\pi_*\beta \ne 0$ in $H_2(X,\ZZ
)$. Then we have:
\begin{enumerate}
\item Every nonempty fiber of $\sigma$ is isomorphic to an open
subset of a projective space.

\item If $c_H(\beta)\ge -1$, then $\sigma$ is a smooth morphism at
general points.

\item If $c_H(\beta)\ge 0$, then the general fiber of $\sigma$ is
nonempty.

\item If $c_H(\beta)\ge 0$ and $\dim Z_i=0$ $\forall\, i$, then $\sigma$ is surjective.

\end{enumerate}}

%

\begin{proof} First note that $\mathrm{Mor}_{\pi _*\beta}
(\mathbb{P} ^1, X)$ contains $\mathrm{Mor}_\beta (\mathbb{P} ^1,
\tilde{X})^\sharp$ as a quasi-projective subvariety over which the
scheme $\mathrm{ev}^{-1}(Z_i)$ has the relative Hilbert polynomial
$e_i$ for all $i=1,...,r$. We will describe a fiber of $\sigma$ as
a subscheme in
\[ \mathrm{Mor}_{\pi _*\beta} (\mathbb{P} ^1, X)\subset
\PP(H^0(\PP ^1, \mathcal{O}_{\mathbb{P} ^1} (d) \otimes \CC
^{n+1}) )\] where $d=(\pi _*\beta ) \cdot H$. If we let
\[\begin{array}{c}
P:=\prod P_i \in \prod \mathrm{Hilb} ^{e_i}\PP ^1\times Z_i,
\ \ \ P_i=\sum _{a} e^{(i)}_a(p^{(i)}_a,q^{(i,a)}),\\
\mathrm{Supp}(P_i) = \{ (p_a^{(i)},q^{(i,a)}) \}_a ,\  p_a^{(i)} \ne
p^{(i')}_{a'} \mbox{ if }(i,a)\ne (i',a'),\  q_0^{(i,a)}\ne 0,
\forall\, (i,a)
\end{array}\] with the condition $e^{(i)}_a =1, \, \forall\, a$ if
$\dim Z_{i'} \ne 0$ for some $i'$, then $\sigma ^{-1}(P)$ is a
subvariety of $\PP H^0(\PP ^1,K _{P})$ where $K_{P}$ is the kernel
of the morphism of sheaves
\[\begin{array}{ccc}\mathcal{O}_{\mathbb{P} ^1} (d) \otimes \CC ^{n+1}& \rightarrow &
(\bigoplus _{i,a}\mathcal{O}_{p_a^{(i)},\PP
^1}/m_{p^{(i)}_{a}}^{e^{(i)}_a} )\otimes \mathcal{O}_{\PP ^1}
(d)\otimes \CC
^{n}\\
f &\mapsto & \sum _{j=1}^n (\oplus _{i,a} [q_0^{(i,a)}f_j -
q_j^{(i,a)} f_0])\otimes 1_j, \end{array}\] where $1_j =
(\overbrace{0,...,1}^{j},0,...,0) \in \CC ^n$. Precisely speaking,
$\sigma ^{-1}(P)_{\mathrm{red}}$ coincides with \[(*)\ \ \ \PP
H^0(\PP ^1,K _{P}) \cap \mathrm{Mor}_{\pi _*\beta} (\mathbb{P} ^1,
X) \setminus \bigcup _{ P'\in \mathrm{Hilb}^{e+1}(\PP ^1\times Z)
\,:\, P'\supset P} \PP H^0 (\PP ^1, K_{P'}),
\] where if $P'$ is not simple at $(p,q)$ and $\dim Z_i\ne 0$ for
some $i$, then $K_{P'}$ is defined as the kernel of
\[ \begin{array}{ccc} K_P &\rightarrow & m_p/m_p^2
\otimes \mathcal{O}_{\PP ^1} (d) \otimes N_{Z/X}|_q \\
f & \mapsto & \sum _j(q_0f_j-q_jf_0)\otimes
[\frac{\partial}{\partial z_j}],
\end{array}\] where $N_{Z/X}|_q:=T_X|_q /\,
T_Z|_q$ (normal space) and $\{ z_j:=x_j/x_0 : j=1,...,n\}$ are the
coordinates of $\CC ^n = \{ x_0\ne 0 \} \subset \PP ^n$.

Since $\mathcal{O}_{\PP ^1}(d- \sum _{i,a} e_a^{(i)}
p_a^{(i)})\otimes \CC ^{n+1} \subset K_P \subset \mathcal{O}_{\PP
^1}(d)\otimes \CC ^{n+1}$, the sheaf $K_P$ is isomorphic to
$\oplus _{j=0}^ n \mathcal{O}_{\PP ^1}(k_j)$ for some $k_j$ with
constraints, $d\ge k_j \ge d-\sum _{i=1}^r e_i$ for all $j$. Now
when $\dim Z_i=0$ for all $i$ and $k_j\ge -1$ for all $j$, then
$\dim \PP H^0(\PP ^1, K_P) = \dim \PP H^0(\PP ^1,\mathcal{O}_{\PP
^1}(d)\otimes \CC ^{n+1}) - ne$. When $\dim Z _i >0$ for some $i$
and $k_j\ge 0$, then $\dim\PP H^0(\PP ^1, K_{P'}) = \dim \PP
H^0(\PP ^1,\mathcal{O}_{\PP ^1}(d)\otimes \CC ^{n+1}) - n(e+1)$
for $P'\supset P$. These facts applied to $(*)$ complete the
proof. \end{proof}

{\Remark\label{product-fiber} Note that Lemma \ref{fiber} holds
true for a product $X=\prod _k \PP ^{n_k}$ of projective spaces if
we let $c_H(\beta ) := \mathrm{min} \{ \beta \cdot (\pi ^*H_k-
E)\}_k$ where $H_k$ is the hyperplane class of $k$-th component of
the product space $X$.}

\subsection{Some elementary facts}

The followings are standard facts.

{\Proposition\label{deform}{\rm(cf. [8, Proposition II.3.7],
\cite{Deb})} Let $X$ be a smooth variety and $Y$ be a subvariety.
Then any free morphism $f:\PP ^1\ra X$ can be deformed to a
morphism $f_{\epsilon}:\PP ^1\ra X$ which is transversal to $Y$.}

{\Lemma\label{irred} Let $X$ and $Y$ be varieties and assume that
$Y$ is irreducible. Let $f:X\rightarrow Y$ be a dominant morphism
in any irreducible component of $X$. Then if the general fiber of
$f$ is irreducible, then $X$ is irreducible.}

\begin{proof} The proof is straightforward. \end{proof}

\subsection{A consequence}

Let $X$ be a product $\prod _k \mathbb{P} ^{n_k}$ of projective
spaces $\PP ^{n_k}$ and let $\tilde{X}$ be a blowing-up of $X$
along a smooth closed subvariety $Z$. Denote by $E$ the
exceptional divisor and denote by $H_k$ the divisor class. We
assume that $\pi _*\beta \ne 0$ and $e_i\ge 0$ for all $i$.

{\Theorem\label{theoremA}
\begin{enumerate} \item  If $\beta \cdot
(\pi ^*H_k - E) \ge -1,$ for all $k$ and $Z$ are finite points,
then $\mathrm{Mor}_\beta (\mathbb{P} ^1, \tilde{X})$ is an
irreducible smooth variety.

\item If $\beta \cdot (\pi ^*H_k - E) \ge 0$ for all $k$, then
$\mathrm{Mor}_\beta (\mathbb{P} ^1, \tilde{X})^\sharp$ is a
nonempty, irreducible smooth variety.

\item If $\beta \cdot (\pi ^*H_k - E) \ge 0$ for all $k$, and all
centers $Z_i$ are convex, then $\mathrm{Mor}_\beta (\mathbb{P} ^1,
\tilde{X})$ is an irreducible smooth variety. \end{enumerate} }

\begin{proof} We prove it when $X=\PP ^n$. The condition $\beta \cdot (\pi
^*H - E) \ge -1$ implies that $H^1(\mathbb{P} ^1, (f^*T_X)(-e))=0$
by the Euler sequence on $\mathbb{P}^n$. Now by (1) of Lemma
\ref{exact}, $H^1(\mathbb{P} ^1, (\tilde{f}^*T_{\tilde{X}}))=0$,
which implies that every irreducible component of
$\mathrm{Mor}_\beta (\PP ^1, \tilde{X})^\sharp$ is smooth with the
expected dimension.

{\em Proof of} (1). The first assertion follows from Lemma
\ref{fiber} and Lemma \ref{irred}.

{\em Proof of} (2). Let us prove the second assertion of the
theorem. Every element in $\mathrm{Mor}_\beta (\PP ^1,
\tilde{X})^\sharp$ is a free morphism by (1) of Lemma \ref{exact}.
Therefore by Proposition \ref{deform}, it is enough to show the
irreducibility of the sublocus $\mathrm{Mor}_\beta (\PP ^1,
\tilde{X})^\circ$ of morphisms which are transversal to $E$. By
Lemma \ref{fiber} the general fiber of $\sigma$ restricted to any
component of $\mathrm{Mor}_\beta (\PP ^1, \tilde{X})^\circ$ has
the expected dimension and irreducible. Thus the morphism $\sigma$
restricted to any component is dominant on the irreducible variety
$\prod \mathrm{Hilb}^{e_i}(\PP ^1\times Z_i)$. Now the proof
follows from Lemma \ref{irred}. The moduli space is nonempty by
Lemma \ref{fiber}.

{\em Proof of} (3). If a morphisms $\tilde{f}$ lie on $E$, using
Lemma \ref{ruled} with $k=0$, we can deform it to an element in
$\mathrm{Mor}_\beta (\PP ^1, \tilde{X})^\sharp$: By Lemma
\ref{ruled}, $H^1(\PP^1, \tilde{f}^*T_{\tilde{X}})=0$ which
implies the moduli space $\mathrm{Mor}_\beta (\PP ^1, \tilde{X})$
is smooth at $\tilde{f}$. By the exact sequence in iii) of the
proof of Lemma \ref{ruled}, we see that there is a deformation of
$\tilde{f}$ to an element in $\mathrm{Mor}_\beta (\PP ^1,
\tilde{X})^\sharp$. Now the proof completes by 2) above.

\medskip

By Remark \ref{product-fiber}, the same proof for the case of the
product of projective spaces works. \end{proof}

\section{Successive blowing-up Case}
\subsection{Set-up}
Let $X_0=X$ be a smooth projective variety and let $\pi _i: X_i
\rightarrow X_{i-1}$, $i=1,...,r$, be a blowing-up of $X_{i-1}$
along a smooth irreducible subvariety $Z_i$. In general, the space
$X_r$ is a successive blowing-up of $X$. Let $E_i^t\subset X_r$
(resp. $E_i^s\subset X_r$) be the total (resp. strict) transform
of the exceptional divisor associated to $Z_i$ and let
$e_i=\beta\cdot E_i^s$. Denote by $\mathrm{Mor}_\beta (\mathbb{P}
^1,X_r)^\circ$ the sublocus of $\mathrm{Mor}_\beta (\mathbb{P}
^1,X_r)$ of the morphisms $\tilde{f}$ which are transversal to
$E=\cup _{i=1}^r E_i^s$ and do not intersect with $E_i^s\cap
E_j^s$ for $i\ne j$. Then for $e_i\ge 0, \ \forall\, i$ we obtain
a morphism
\[\begin{array}{ccc} \sigma: \mathrm{Mor}_\beta (\mathbb{P} ^1,
X_r)^\circ &\rightarrow& \prod
_{i=1}^r\mathrm{Hilb}^{e_i}(\mathbb{P}
^1\times Z_i) \\
\tilde{f} &\mapsto& \prod _{i=1}^{r} (\Gamma _{\pi _{i}\circ
...\circ \pi _r \circ \tilde{f}}\cap \PP ^1\times (Z_i\setminus
\cup _{j>i}(\pi _i\circ ...\circ \pi _r (E_j^s))
\end{array}
\] as the generalization of the previous $\sigma$ in subsection
\ref{sigmasubsection}.

Inductive application of the exact sequence of Lemma \ref{exact}
proves the following corollary.

{\Corollary\label{smooth} Suppose that $H^1(\mathbb{P} ^1,
f^*(T_X(-\sum E_i^t)))=0$, where $f=\pi _1\circ ...\circ \pi
_r\circ\tilde{f}$. Then the morphism $\sigma$ is smooth at
$\tilde{f}$.}

\subsection{Jet spaces}\label{jetspaces}

We want to show that the general fiber of $\sigma$ is rationally
connected, provided with a suitable condition on $\beta$ when
$X=\PP ^n$ or their products. However it is hard to analyze the
fiber of $\sigma$ directly as done in subsection
\ref{fibersubsection}. Our strategy is to introduce an auxiliary
morphism $\tau$ by imposing further conditions on jets of the
morphism $f: \mathbb{P}^1 \rightarrow X$. It turns out that the
fiber of $\tau$ is simple to study. Since the jet conditions on
$f$ can be translated to the vanishing conditions on the
blowing-up space, we will be able to express the general fiber of
$\sigma$ by the fibers of $\tau $ (more precisely its product
$\tau _{\bf m}$) which are rationally connected.

To introduce $\tau$, let $J^k_q X=\mathrm{Mor}((\mathrm{Spec}\CC
[\epsilon ]/(\epsilon ^{k+1}),0),(X,q))$ be the $k$-jet space of
$X$ at $q \in X$. Then the morphism $f:\mathbb{P} ^1 \rightarrow
X$ naturally assigns an element $[f]^k_p \in J^k_{f(p)}X$ for any
$p\in \mathbb{P} ^1$. Using the assignment we define a morphism
\[\begin{array}{ccc} \tau: ((\mathbb{P} ^1)^l \setminus
\Delta )\times \mathrm{Mor}_\beta (\mathbb{P} ^1,X) &\rightarrow&
(\mathbb{P} ^1 \times J^k X)^l \\
({\bf p},f) &\mapsto& (p_i, [f]^k_{p_i})_{i=1,...,l},
\end{array}\]
where $\Delta$ is the big diagonal and $J^kX = \coprod _{p\in
X}J^k_p X$.

{\Lemma\label{smoothnessoftau} If $H^1(\mathbb{P} ^1, f^*T_X
(-(k+1)l))=0 $, then $\tau$ is smooth at $({\bf p},f)$.}

\begin{proof} The natural exact sequence $$0\rightarrow f^*T_X
(-(k+1)l)) \rightarrow f^*T_X \rightarrow f^*T_X \otimes
(\bigoplus _{i=1}^l \mathcal{O}_{p_i,\PP ^1}/m_{p_i}^{k+1})
\rightarrow 0$$ induces the map $$ H^0(\PP ^1, f^*T_X) \rightarrow
H^0(\PP ^1, f^*T_X \otimes (\bigoplus _{i=1}^l
\mathcal{O}_{p_i,\PP ^1}/m_{p_i}^{k+1}))$$ which is the tangent
map $T\tau _{(\prod p_i,f)}|_{0 \times H^0(\PP ^1, f^*T_X)}$.
(Note that the exactness holds since $p_i$ are pairwise distinct
for all $i$.) Indeed
$$T_{J^k X}|_{[f]^k_p} = H^0 (\Spec \CC [\epsilon ]/(\epsilon
)^{k+1}, ([f]^k_p)^* T_X)) = H^0 (\PP ^1, f^*T_X\otimes
\mathcal{O}_p/m_p ^{k+1})$$ by base change. This induced morphism
is surjective by the assumption. \end{proof}

\subsection{A general simple fact}
Let $\pi :\tilde{X} \rightarrow X$ be a blowing-up of a smooth
variety $X$ along a subvariety $Z$. Let $E_q:= \pi ^{-1}(q)$ where
$q$ is a smooth point of $Z$. Note that there is a natural
morphism $j^k_q: (J_q^kX)^\circ \ \rightarrow\bigcup _{w\in E_q}
J_w^{k-1}\tilde{X}$ defined by lifting, where $(J_q^kX)^\circ$
consists of $s \in J^k_qX$ which are, as $k$-jet arc, transversal
to $Z$ at $q$.

{\Lemma\label{jet} The morphism $j^k_q$ is smooth and every fiber
is a rational variety. }

\begin{proof} This is a local problem at $q$. So we may assume that
$X=\mathbb{C} ^n$, $q=0$, $Z=\{ 0\}\times \mathbb{C} ^l$ and $\pi
(t,{\bf x},{\bf y}) =(t,t{\bf x},{\bf y})$. Consider a $k$-jet
$s(t)$ at $t=0$ such that $\frac{d s_1}{dt}|_{t=0}\ne 0$, then
$(s_1(t),\frac{s_2(t)}{s_1 (t)},...,\frac{s_{n-l}(t)}{s_1 (t)},
s_{n-l+1}(t),...,s_n(t))$ mod $t^k$ is, by definition,
$j^k(s(t))$. This shows that the morphism $j^k$ is regular. Then
it is straightforward to check the smoothness of $j^k_q$ and the
rationality of the fiber. \end{proof}

\subsection{The morphism $\tau _{\bf m}$ and its fibers when $X=\prod _{j=1}^m\mathbb{P}
^{n_j}$}

We define a morphism
$$\tau _{\bf m}: ((\mathbb{P} ^1)^{\sum e_i}
\setminus \Delta )\times \mathrm{Mor}_{\beta} (\mathbb{P} ^1,X)
\rightarrow \prod _{i=1}^{r} (\mathbb{P} ^1\times J^{m_i}X )^{e_i}
$$ similar to the  $\tau$ as in subsection \ref{jetspaces}. Here $m_i$ are nonnegative integers and ${\bf m}=(m_1,...,m_r)$.

{\Lemma\label{rationalityoftau} If the target space $X$ is a
product of projective spaces, $\prod _{j=1}^m\mathbb{P} ^{n_j}$,
then the fibers of $\tau _{\bf m}$ (with their induced reduced
scheme structure) are rational varieties.}

\begin{proof} We will prove that when $X=\mathbb{P} ^n$ and $l=1$, every
fiber of $\tau$ defined subsection \ref{jetspaces} is a linear
subvariety of $\PP (H^0(\PP ^1,\mathcal{O}(d))\otimes \CC
^{n+1})$. The general case is straightforward from the proof. Now
it is easy to check that $p\times f$ and $p\times g$ are in a same
fiber of $\tau$, then $p\times (\mu f + \lambda g := (\mu f_0
+\lambda g_0 ,..., \mu f_n+\lambda g_n))$ is in the same fiber for
all but finite number of $(\mu ,\lambda )\in \mathbb{P} ^1$.
\end{proof}

\subsection{A proof of the Main Theorem and example}

Let $m_i=\# \{Z_j \ | \ j< i ,\ \ (\pi _j\circ ...\circ \pi
_r)^{-1}(Z_j)\supset E_i^t\}$. So general points of $Z_i$ are the
$m_i$-th infinitesimal points of $X$. For each $i$, we reindex
$Z_j$ so that  $\pi ^{-1} (Z_{i_k})\supset E_i$ where $i_1 < ... <
i_{m_i}$.

\bigskip

\noindent{\em Proof of {\rm (1)} and {\rm (2)} of the Main
Theorem.} By the assumption on $\beta$ and (1) of Lemma
\ref{exact}, the space $\mathrm{Mor}_\beta (\mathbb{P} ^1,
X_r)^\sharp$ is smooth and has the expected dimension, and its
elements are free morphisms. Hence $\mathrm{Mor}_\beta (\mathbb{P}
^1,X_r)^\circ$ is open dense in $\mathrm{Mor}_\beta (\mathbb{P}
^1,X_r)^\sharp$ by Proposition \ref{deform}. Also note that
$\sigma$ is a smooth morphism by Corollary \ref{smooth}. Therefore
by Lemma \ref{irred} and Theorem \ref{GHS} below, in order to
verify (1) and (2) of the Main Theorem, it is enough to show that
the general fiber of $\sigma $ is rationally connected. Let
$P=\sum (p_i,q_i) \in \prod \mathrm{Hilb}^{e_i}(\mathbb{P}
^1\times Z_i)$ such that all points in $P$ are simple. Then we
obtain the inclusion
$$(**) \ \ \ \pi _\beta \circ \sigma  ^{-1}(P)
\subset \mathrm{pr}_2 \circ(\tau _{\bf m}) ^{-1} (( p_i \times
(j^1_{\pi_{i_{m_i}}(q_i)}\circ ...\circ j^{m_i}_{\pi _{i_1}\circ
...\circ \pi _{i_{m_i}}(q_i)} )^{-1} (q_i)) _{i=1}^{\sum e_j})$$
where $\pi _\beta : \mathrm{Mor}_\beta (\mathbb{P} ^1, X_r)^\circ
\subset \mathrm{Mor}_{\pi _*\beta} (\mathbb{P} ^1, X)$ is the
natural embedding, and $\mathrm{pr} _2$ is the projection to the
second factor. In $(**)$ the right hand side (for short, RHS)
includes morphisms $\tilde{f}$ with $\deg \tilde{f}^{-1}(E)\ge
\beta\cdot E$. This is the reason why both sides may not coincide.

Since $\tau _{\bf m}$ is a dominant morphism with a rationally
connected general fiber by Lemma \ref{smoothnessoftau} and Lemma
\ref{rationalityoftau}, RHS is also rationally connected, thanks
to Theorem \ref{GHS}. Since LHS is open subset of RHS, we conclude
that LHS is also rationally connected.

\noindent{\em Proof of {\rm (3)} of the Main Theorem.} In case i):
If $f$ lies on $E_i^s$ for some $i$, then $e_i<0$: Take a
hypersurface $Y$ of $X$ such that $Y$ contains the image of $f$;
and the strict transform $D$ of $Y$ under $\pi$ does not contain
the image of $\tilde{f}$. Then $0=Y\cdot f_*[\PP ^1 ]=\pi ^* Y
\cdot \tilde{f}_*[\PP ^1] = (D + \sum a_i E_i^s )\cdot
\tilde{f}_*[\PP ^1]$ shows that $E_i^s\cdot \tilde{f}_*[\PP ^1]$
is negative for some $i$.

In case ii): This is the case (3) of Theorem \ref{theoremA}.
\hfill$\Box$

{\Theorem {\rm (\cite{GHS})}\label{GHS} Let $f: X \rightarrow Y$
be a dominant morphism between irreducible varieties $X$ and $Y$.
If \/ $Y$ and the general fiber of $f$ is rationally connected,
then $X$ is rationally connected.}

\bigskip

{\Example Let $Y$ be a quadratic line complex in $\mathbb{P} ^5$,
which is a complete intersection of two smooth quadrics in
$\mathbb{P}^5$. Then $Y$ is isomorphic to the moduli space of
isomorphism classes of stable rank 2 vector bundle on a curve of
genus $g=2$ with fixed determinant of degree 1 \cite{New}. Let
$\tilde{X}$ be the blowing-up of $Y$ along a line and let $\tilde
Q$ be the inverse image of the line. Then $\tilde{X}$ is a
blowing-up of $\mathbb{P} ^3$ along a smooth quintic curve $C$.
Let $E_C$ be the inverse image of the curve $C$.
$$\begin{array}{ccccc}
& E_C\subset & \tilde X & \supset\tilde Q & \\
& \pi \swarrow & & \searrow \pi ' & \\
C\subset\mathbb{P} ^3 & & & & Y\supset\mbox{line}\\
\end{array}$$
Then $\pi (\tilde Q)$ is a quadric surface $Q$ containing $C$.
Therefore $\tilde Q=2H-E_C$ where $H$ is the proper transform of a
hyperplane class in $\mathbb{P}^3$. Castravet [Cas] shows that
there are at least two components (nice one, and almost nice one)
of $\mathrm{Mor}_d(\mathbb{P} ^1, Y)$ with the expected dimension.
The almost nice component consists of morphisms $\mathbb{P}^1\to
Y$ which are $d$ to 1 onto lines in $Y$ \cite{Cas}. Here $0< d\in
\ZZ \cong H_2(Y,\ZZ )$ with respect to the ample generator of
$\mathrm{Pic}Y$. These two components $\mathrm{Mor}_d(\mathbb{P}
^1, Y)$ are birational to two components of $\mathrm{Mor}_{(d,e)}
(\mathbb{P} ^1, \tilde{X})$ with $e=2d$. Here $(d,e)\in \ZZ \times
\ZZ \cong H_2(\tilde{X},\ZZ)$ with respect to $\pi ^*(H)$ and
$E_C$. In this case note that at every point in the corresponding
component of the almost nice component, $\sigma$ is not smooth.}

\end{document}